\documentclass[11pt]{article}
\usepackage{float,amsmath,amsthm,amssymb,mathrsfs,enumerate,graphics,tikz-cd,geometry,todonotes}
\usepackage{hyperref}

\usetikzlibrary{tikzmark,calc}
\usepackage{algorithm}
\usepackage[noend]{algpseudocode}

\algrenewcommand\algorithmicwhile{\textbf{While}}
\algrenewcommand\algorithmicfor{\textbf{For}}
\algrenewcommand\algorithmicdo{\textbf{Do}}
\algrenewcommand\algorithmicif{\textbf{If}}
\algrenewcommand\algorithmicthen{\textbf{Then}}
\algrenewcommand\algorithmicelse{\textbf{Else}}
\algrenewcommand\algorithmicend{\textbf{End}}
\algrenewcommand\algorithmicreturn{\textbf{Return}}

\newcommand{\QQ}{\mathbb{Q}}

\newcommand{\PP}{\mathbb{P}}

\newtheorem{theorem}{Theorem}[section]
\newtheorem{lemma}[theorem]{Lemma} 
\newtheorem{prop}[theorem]{Proposition}

\newtheorem{remark}[theorem]{Remark}  
\newtheorem{example}[theorem]{Example}

\title{How to Reconstruct a Planar Map From its Branching Curve}
\author{Eriola Hoxhaj and Josef Schicho \\ RISC, Johannes Kepler University of Linz, Austria }
\date{}

\begin{document}

\maketitle

\begin{abstract}
We present an algorithm for constructing a map $\PP^2\to\PP^2$ with a given branching curve.
The stepping stone is the ramification curve, which is obtained as the linear normalization of the branching curve.
\end{abstract}

\section*{Introduction}
We provide an algorithm that deals with the following reconstruction problem: given a plane curve $B \subset \PP^2$,
compute a map $f :\PP^2 \to \PP^2$ such that $B$ is the branching curve of this map, assuming that the map is generic.
By \cite{Ciliberto_Flamini:11}, genericity implies that the ramification curve is smooth, the branching curve has only nodes
and cusps, and the map from the ramification curve to the branching curve is birational.

The problem is motivated by potential applications in computer vision: a picture of a parametrised surface can be mathematically
described by the composition of the parameterization followed by a central projection to the image plane, and the apparent contour
of the parameterised surface is the branching curve of the composed map. A solution to the reconstruction problem, if it is unique,
should allow to obtain information about the parameterised surface.

The uniqueness of the solution, up to unavoidable projective transformations, is guaranteed by a much more general result in
\cite{Kulikov:99}: any generic map $f:S_1\to S_2$ of degree at least $5$ is determined up to right composition with
isomorphism to $S_1$ by its branching curve in $S_2$. In our case, the map $f:\PP^2\to\PP^2$ is defined by polynomials of degree~$d$
that do not have a common zero, and this implies that the degree of the map is $d^2$. For $d\ge 3$, uniqueness is therefore
guaranteed; our constructive algorithm could be considered as another proof for this special case. The case $d=2$ is classical.
Here the branching curve is a sextic with 9 cusps, dual to a smooth cubic. Chisini~\cite{Chisini:44} showed that there are four maps
$S_1\to\PP^2$, up to right equivalence. In three of these four solutions, $S_1$ is isomorphic to $\PP^2$. The three solutions turned
out to be left-right-equivalent, i.e., there are isomorphic up to composing with an isomorphism on the right and with an automorphism
on the left that preserves the branching curve.

Our algorithm computes first the linear normalization $R$ of the given branching curve (see \cite{GKZ} for a discussion of the
linear normalization of a projective variety). We show that in our case -- assuming that the reconstruction problem has a solution --
the linear normalization is the ramification curve of a generic projection from the Veronese surface $V_d\in\PP^{\frac{(d(d+3)}2}$
to $\PP^2$. The computation of the linear normalization can in theory be achieved by computing the normalization of the cone over
the branching curve. However, the curve $B$ has degree $3d(d-1)$ and $\frac{9(d^2-2)(d-1)^2}{2}$ singularities (nodes and cusps),
so the algorithms available in special purpose computer algebra systems would not be fast enough even for the case $d=3$. Hence we provide a new algorithm
for the linear normalization, adapted to the situation of curves that have only nodes and cusps. It works well for $d=3$ and $d=4$.

Once we have the ramification curve $R\in \PP^{\frac{d(d+3)}2}$ (in case $d\ge 3$), we show that there is a unique Veronese surface $V_d$
containing $R$, and we compute the equations of this Veronese surface. This step uses syzygies in a way similar to \cite{Bothmer:2008},
where a curve is given too, and certain specific varieties containing the curve need to be constructed. The classical case $d=2$ is treated
by a different approach, using the Hessian pencil of the dual curve (see \cite{Arbetani_Dolgachev:09}).

The Veronese surface $V_d$ is isomorphic to $\PP^2$. The last step of our algorithm is to compute an explicit isomorphism $\PP^2\to V_d$,
i.e., a parametrisation of $V_d$. This is even more expensive than the computation of the linear normalization, see for instance \cite{Schicho:05c} for the case $d=3$. The hardest part
is a Diophantine subproblem: one needs to find a point in $V_d$ first. Unfortunately, the available algorithms for solving this Diophantine problems over $\QQ$ are too slow to handle random examples with $d=3$. Hence this substep of the algorithm is just theoretical. Here we assume that we can find two distinct rational points and proceed to computing a parametrization. We illustrate this step by examples in positive characteristics: here it is easy to find rational points probabilistically, by intersecting with random linear subspaces of codimension two. 

Here is an outline of the paper. Section~1 gives a precise specification of the problem, including the equivalence relation we need to take into account in order to speak of uniqueness of the solution. Section~2 is dedicated to the case $d=2$. In Section~3, we give our new algorithm
for linear normalization and use it to construct the ramification curve $R$. Section~4 is dedicated to the construction of the unique
Veronese surface containing $R$. In Section~5, we give a parametrization algorithm and put the other parts together.

\section{Generic Planar Maps and their Branching Curves}

Let $k$ be a computable field.
A regular map  $f:\PP^2\to \PP^2$ is a called a {\em planar map}. Every planar map is defined by three homogeneous polynomials
$F_0,F_1,F_2\in k[s,t,u]$
of the same degree $d$ without a common zero; we say that $d$ is the {\em polynomial degree} of the map $f$. This number is not
equal to the degree of $f$, i.e., the cardinality of the preimage of a generic point. Specifically, the degree of the map is $d^2$,
by Bezout's theorem.

Let $f:\PP^2\to \PP^2$ be a planar map of polynomial degree $d\ge 2$. The {\em ramification curve} $R$ of $f$ is the scheme defined by
the homogeneous Jacobi matrix $\det\left(\frac{\partial{(F_0,F_1,F_2)}}{\partial{(s,t,u)}}\right)$. This is a curve of degree $3(d-1)$ (maybe reducible,
maybe with multiple components). The {\em branching curve} $B$ is defined as the image of the ramification curve. We say that the planar map $f$
is {\em generic} if it satisfies the following conditions:

\begin{itemize}
\item the ramification curve $R$ is nonsingular (in particular irreducible);
\item the restriction of $f$ to $R$ is birational to $B$;
\item the branching curve $B$ has only nodes and cusps.
\end{itemize}

It is well-known that the set of coefficients vectors of planar maps that satisfy the conditions above is a Zariski-open and dense subset in
the set of all coefficient vectors (see \cite{Ciliberto_Flamini:11}). From now on, assume that $f$ is a generic planar map.
Then $\deg(B)=d\deg(R)=3d(d-1)$.
Since $R$ is smooth, its geometric genus is equal to the arithmetic genus $\frac{(3(d-1)-1)(3(d-1)-2)}{2}$.
The geometric genus is a birationally invariant, thus the geometric genus of $B$ is the same as the geometric genus of $R$.
From the genus formula for planar curves, we can find the number of singularities (in our case nodes and cusps)
\[ \frac{(3d(d-1)-1)(3d(d-1)-2)}{2}-\frac{(3(d-1)-1)(3(d-2)-2)}{2}=\frac{9(d^2-2)(d-1)^2}{2} . \]

We say that two planar maps $f_1,f_2:\PP^2\to\PP^2$ are {\em right equivalent} if there is a projective isomorphism $t:\PP^2\to\PP^2$
such that $f_2=f_1\circ t$.
Similarly, we say that  $f_1,f_2:\PP^2\to\PP^2$ are {\em left equivalent} if there is a projective isomorphism $t:\PP^2\to\PP^2$
such that $f_2=t\circ f_1$.
Left equivalent maps have the same ramification curve and projectively equivalent branching curves.
Right equivalent maps have the same branching curve and projectively equivalent ramification curves.

The goal of this paper is to provide an algorithm for the following problem:

\floatname{algorithm}{Problem}
\renewcommand{\thealgorithm}{DiscToMap} 
\begin{algorithm}[H]
\caption{}\label{algorithm:problem}
\begin{algorithmic}[1]
  \Require an irreducible curve $B$ of degree $3d(d-1)$ in $\PP^2$ that has only nodal and cuspidal singularities.
  \Ensure a generic planar map $f:\PP^2\to\PP^2$ that has $B$ as a branching curve, if such a map exists.
\end{algorithmic}
\end{algorithm}

If the no generic planar map exists such that the branching curve is the given curve, then the algorithm detects this. We also prove that
for $d\ge 3$, the solution is unique up to right equivalence; but, as mentioned in the introduction, the uniqueness is already
a consequence of Kulikov's result~\cite{Kulikov:99}.

The first step in our algorithm will be to compute the {\em linear normalization} of $B$, which we now define,
following~\cite{GKZ}.
For two integers $n>m>2$, a rational map $\PP^n\dashrightarrow\PP^m$ given by an $(m+1)\times(n+1)$ matrix of full rank is called a projection;
its center is the projectivised nullset of the defining matrix.
Let $X\subset\PP^n$ be an irreducible variety that is not contained in a hyperplane.
The projection $p:\PP^n\dashrightarrow\PP^m$ is called a {\em regular projection} if the center of $p$ does not intersect $X$
and $p$ maps $X$ birationally to its image $p(X)$. Regular projections preserve the degree, i.e., $\deg(X)=\deg(p(X))$.
The {\em linear normalization} of $X$ is an irreducible variety $\tilde{X}$ in some $\PP^N$ together with a regular projection
$p:\PP^N\dashrightarrow\PP^n$ that maps $\tilde{X}$ to $X$ and which is universal in the following sense: for any irreducible variety $X'$
in some $\PP^M$ together with a regular projection $p':\PP^M\dashrightarrow\PP^n$ that maps $X'$ to $X$, there is a unique
regular projection $q:\PP^N\dashrightarrow\PP^M$ that maps $\tilde{X}$ to $X'$ and that satisfies $p=p'\circ q$.
It is well-known that the linear normalizations always exists and is unique up to right equivalence.
If the linear normalization is a projective isomorphism, then we say that $X$ us linearly normal.
A normal variety $X\subset\PP^n$ is linearly normal if and only if the restriction map of linear forms
$\Gamma(\PP^n,{\cal O}_{\PP^n}(1))\to\Gamma(X,{\cal O}_X(1))$ is surjective.

\begin{example} \label{ex:nc} \rm
The linear normalization of the nodal cubic in $\PP^2$ with equation $y^2z-x^2z-x^3=0$ is the twisted cubic in $\PP^3$
with equations $zy-wx=w^2-z^2-xz=wy-xz-x^2=0$. The regular projection is $(x:y:z:w)\mapsto(x:y:z)$ with center $(0:0:0:1)$.
\end{example}

The linear normalization is useful because it is projectively isomorphic  to the ramification curve; the precise statement is below.
We write $\nu_d:\PP^2\to\PP^{\frac{d(d+3)}{2}}$ for a Veronese map defined by a basis of forms if degree $d$, and $V_d$ for the Veronese surface
$\nu_d(\PP^2)$. The map $\nu_d$ is unique up to left equivalence.

\begin{prop} \label{prop:veronese}
Assume that $B$ is the branching curve of a generic planar map $f:\PP^2\to\PP^2$ of polynomial degree $d\ge 2$.
Let $\tilde{B}\subset\PP^N$, $p:\PP^N\dashrightarrow\PP^2$ be the linear normalization of $B$.
Then $N=\frac{d(d+3)}{2}$, and there is a Veronese embedding $\nu_d:\PP^2\to V_d\subset\PP^N$ such that the diagram
$$\begin{tikzcd}
 & \tilde{B}\subset V_d\subset\PP^N \arrow[dr,"p"] \\
R\subset\PP^2\arrow[ur,"\nu_d"] \arrow[rr,"f"] && B\subset \PP^2
\end{tikzcd}
$$
commutes. The Veronese embedding $\nu_d$ maps the ramification curve $R$ of $f$ isomorphically to $\tilde{B}$. In particular,
$\tilde{B}$ is contained in $V_d$.
\end{prop}

\begin{proof}
Let $B'\subset\PP^\frac{d(d+3)}{2}$ be the image of $R$ under the Veronese map. The map $f$ is defined by three forms of degree $d$,
hence there is a projection map $p':\PP^\frac{d(d+3)}{2}\to\PP^2$ such that $f=p'\circ \nu_d$. The center of the projection $p'$ does
not intersect $V_d$, because the intersection is the image of the base locus of $f$ and this base locus is empty by assumption.
Hence $p'$ is a regular projection from $B'$ to $B$. By the universality property of the linear normalization, there is
a regular projection $q:\PP^N\to\PP^\frac{d(d+3)}{2}$ from $\tilde{B}$ to $B'$.

To show that $q$ is an isomorphism, we show that $R'$ is linearly normal. The restriction map
$\Gamma(\PP^\frac{d(d+3)}{2},{\cal O}_{\PP^\frac{d(d+3)}{2}}(1))\to\Gamma(B',{\cal O}_{B'}(1))$ of linear forms
corresponds to the restriction map $\Gamma(\PP^2,{\cal O}_{\PP^2}(1))\to\Gamma(B,{\cal O}_{B}(d))$
of forms of degree $d$.
This map is surjective because $R$ is a hypersurface and therefore arithmetically Cohen-Macaulay.
\end{proof}

Proposition~\ref{prop:veronese} suggests the following strategy for solving our problem:
\floatname{algorithm}{Strategy for Algorithm}
\renewcommand{\thealgorithm}{DiscToMap} 
\begin{algorithm}[H]
\caption{}\label{alg:main}
\begin{algorithmic}[1]
  \Require an irreducible curve $B$ of degree $3d(d-1)$ in $\PP^2$.
  \Ensure a generic planar map $f:\PP^2\to\PP^2$ that has $B$ as a branching curve.
  \Statex
  \State {\bf Compute} the linear normalization $p:\tilde{B}\to B$ (see section~3).
  \State {\bf Compute} a Veronese surface $V_d$ containing $\tilde{B}$ (see section~4).
  \State {\bf Compute} a Veronese map $\nu_d:\PP^2\to V_d$ (see section~5).
  \State \Return $p\circ \nu_d:\PP^2\to\PP^2$.
\end{algorithmic}
\end{algorithm}
We will follow this strategy for $d\ge 3$. For $d=2$, we could not find a way to do step~3. Fortunately, there
is an alternative, to be described in the next section.

\section{Maps of Polynomial Degree Two}

In this section, we assume that $k$ has characteristic zero and is algebraically closed.

The branching curve $B$ of a generic planar map of polynomial degree two is a sextic of genus one with nine cusps. By the Pl\"ucker formula for
plane curves, the dual curve $\check{B}$ is a nonsingular cubic. On the other hand, the ramification curve $R$ is also a smooth cubic.
Both $R$ and $\check{B}$ are birationally equivalent to $B$, hence they are birationally equivalent to each other. Two cubic curves are
birationally equivalent if and only if they are projectively equivalent; so, we already know $R$ up to projective equivalence.

\begin{lemma} \label{lem:ram}
Every nonsingular cubic is the ramification curve of some generic planar map of polynomial degree two.
\end{lemma}

\begin{proof}
Let $F\in k[s,t,u]$ be a cubic form defining a nonsingular cubic curve. The gradient $\nabla(F)$ defines a planar map of polynomial
degree two. Its ramification curve coincides with the Hessian $H(F)$. Hence the statement is a consequent of the next lemma.
\end{proof}

\begin{lemma} \label{lem:hessian}
Every nonsingular cubic is the Hessian of exactly three other nonsingular cubics.
\end{lemma}

\begin{proof}
We are going to construct a nonsingular cubic form $G$ whose Hessian is a given cubic form $F$ defining a nonsingular cubic curve.
We claim that $G$ can be chosen of the form $aF+bH(F)$ with suitable constants $a,b\in k$.
In order to prove the claim, we assume that $F$ is in Hessian normal form
\[ F_\lambda = s^3+t^3+u^3+\lambda stu, \ \lambda^3+27\ne 0. \]
This is without loss of generality: any nonsingular cubic can be transformed to Hessian normal form by some projective transformation
(see \cite{Arbetani_Dolgachev:09}). The condition $\lambda^3+27\ne 0$ is equivalent to the statement that $\nabla(F)$ does not have zeroes,
or in other words that $F$ defines a nonsingular curve.

The Hessian of a Hessian normal form is also a Hessian normal form: specifically,
\[ H(F_\mu)=F_{-\frac{\mu^3+108}{3\mu^2}}. \]
So we just need to find $\mu$ such that $-\frac{\mu^3+108}{3\mu^2}=\lambda$, and this is always possible.
The fact that the discriminant of $\mu^3+3\lambda \mu^2+108$ with respect to $\mu$ is equal to $\lambda^3+27$ (up to some multiplicative constant)
shows that there are exactly three solutions.
The fact that the resultant of $\mu^3+3\lambda \mu^2+108$ and $\mu^3+27$ with respect to $\mu$ is equal to $\lambda^3+27$ shows that
all three solutions are equations of nonsingular cubic curves.

Assume, indirectly, that $G'$ is the equation of a nonsingular cubic which is not in Hessian normal form such that $H(G')=F$. Let $\Gamma$
be the pencil of cubics generated by $G'$ and $H(G')=F$. The Hessian of any cubic in $\Gamma$ is again in $\Gamma$. In particular, $H(F)$
is in $\Gamma$. Since $F$ and $H(F)$ are distinct, they generate $\Gamma$, and this shows that $\Gamma$ is the pencil of Hessian normal forms.
But this contradicts the assumption that $G'$ is not in Hessian normal form.
\end{proof}

The three solutions are classical~\cite{Chisini:44} -- see also \cite{Catanese:86} for a very explicit discussion including also the fourth map from a surface -- not isomorphic to $\PP^2$ -- to $\PP^2$ with the same branching curve.

Here is the algorithm for the reconstruction of a planar map of polynomial degree~2 that corresponds
to the constructive proof above.

\floatname{algorithm}{Algorithm }
\renewcommand{\thealgorithm}{DiscToMapCase2} 
\begin{algorithm}[H]
\caption{}\label{alg:case2}
\begin{algorithmic}[1]
  \Require an irreducible sextic $B$  with nine cusps.
  \Ensure a generic planar map $f:\PP^2\to\PP^2$ that has $B$ as a branching curve.
  \Statex
  \State {\bf Compute} the dual $\check{B}$ of $B$ (a nonsingular cubic).
  \State {\bf Compute} the Hessian $H(\check{B})$.
  \State {\bf Compute} $\mu\in k$ such that the Hessian of $H(\check{B})+\mu\check{B}$ is equal to $\check{B}$
    (three solutions).
  \State \Return $\nabla(H(\check{B})+\mu\check{B})$.
\end{algorithmic}
\end{algorithm}

\begin{figure}[H]\label{fig:example1}
\centering
\includegraphics[width=6cm,height=5cm]{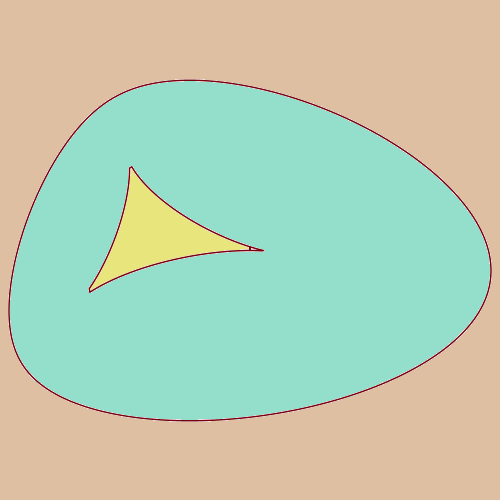}
\caption{The branching curve $B$ separates the plane into three regions. The points in the beige region have no real preimage points for the two first maps and four real preimages for the third map. The points in the blue sky region have two real preimage points and the points in the yellow region have four real preimage points, for all three maps.}
\end{figure}

\begin{example} \rm 
Let $B$ be the curve with equation $$ 32x^6+196x^4y^2+72x^3y^3+436x^2y^4-280xy^5+144y^6$$$$+16x^5z+96x^4yz+352x^3y^2z+140x^2y^3z+700xy^4z-272y^5z$$$$-32x^4z^2+44x^3yz^2+52x^2y^2z^2+184xy^3z^2+209y^4z^2-16x^3z^3$$$$-80x^2yz^3-266xy^2z^3+128y^3z^3+x^2z^4-36xyz^4-132y^2z^4+2yz^5 = 0 .$$
This is a sextic with nine cusps, three of them being real (see Figure~\ref{fig:example1}). Its dual is the cubic curve $\check{B}:t^2u-(s^2-u^2)(s-2u)+ts^2/2=0$. The three numbers such that such that the hessian of $H(\check{B})+\mu\check{B}$ is equal to the $\check{B}$ are $\mu_1=-31.37404937, \mu_2=-14.20865677,\mu_3=45.58270614$.  
We obtain three planar map of polynomials of degree 2 with the sextic $B$ as its branching curve. The picture Figure~\ref{fig:example1} shows the number of real preimages for each connected component of the complement of $B$.
\end{example}

\section{Linear Normalization}

As mentioned in Section~1, 
a normal variety $X\subset\PP^n$ is linearly normal if and only if the restriction map of linear forms
$\Gamma(\PP^n,{\cal O}_{\PP^n}(1))\to\Gamma(X,{\cal O}_X(1))$ is surjective. This suggests an algorithm for the computation
of the linear normalization of $X$: first, compute its projective normalization $N(X)$. Second, take a basis for the linear space
of degree one elements of the graded coordinate ring of $N(X)$, and use this basis to define a map into projective space. The image
of this map is the linear normalization.

This algorithm would be correct, but unfortunately it is too costly for the computation of the linear normalization of branching curves.
Already for $d=3$, the branching curve is a plane curve of degree~18 with 126 singular points with coordinates in a large field extension,
and this is not feasible for the available programs for the computation of the normalization in Singular or Macaulay2. Here we use a more
efficient method that takes advantage of the special situation that computes the degree~1 elements without computing the whole integral
closure. Also, we make use of the assumption that the given plane curve has only nodes and cusps as its singular points. 

Let $k$ be a field and $B\in\PP_k^2$ be an irreducible plane curve defined by a homogeneous polynomial $F\in k[x,y,z]$.
The graded coordinate ring of $C$ is denoted by $S:=k[x,y,z]/\langle F\rangle$. Its integral closure in its quotient field $\mathrm{QF}(S)$ is denoted by $\tilde{S}$.
The ring $\tilde{S}$ is also the graded coordinate ring of a projective curve $N(B)$, the projective normalization of $B$.
The ring inclusion $S\hookrightarrow\tilde{S}$ induces a projection map $N(B)\to B$.

\begin{remark}
It may happen that the graded ring $\tilde{S}$ is not generated by linear elements. Then $N(B)$ is naturally embedded in a weighted
projective space, where the weights are the degrees of ring generators. This paper only contains examples where $\tilde{S}$ is generated in degree~1.
\end{remark}

The conductor ideal $C$ is defined as the conductor of $S$ in $\tilde{S}$, i.e., the set of all elements $s\in S$ such that $s\tilde{S}\subset S$.
It is an ideal in both rings $S$ and $\tilde{S}$. The preimage of $C$ under the quotient map $k[x,y,z]\to S$ is also known as the Gorenstein
adjoint ideal. It is equal to the ideal of all polynomials that vanish with order at least $r-1$ at all $r$-fold singular points of $B$, including
infinitely near singularities (see \cite{JB}). 

\begin{remark} \label{rem:nc}
In our situation, $B$ has only nodes and cups. Hence all multiplicities of points of singular points are equal to $2$, and there are no
infinitely near singularities. It follows that the Gorenstein adjoint ideal is equal to the radical of the singular locus of $B$. 
\end{remark}

\begin{lemma} \label{lem:ch}
Let $a\in\mathrm{QF}(S)^\ast$, and assume that $A\subset\mathrm{QF}(S)$ is a finitely generated $S$-submodule of the quotient field.
If $aA\subseteq A$, then $a$ is integral over $S$.
\end{lemma}

\begin{proof}
Because $A$ is finitely generated, there exists a positive integer $r$ and a surjective map $\pi:S^r\to A$.
Let $\mu:A\to A$ be the multiplication map $m\mapsto am$. Because $\mu$ is $S$-linear and $S^r$ is a projective module, $\mu$ can be lifted
to obtain a commutative diagram
\[ \begin{matrix}
    S^r&\xrightarrow{\quad \psi\quad}& S^r\\
    \bigg\downarrow{\pi}& &\bigg\downarrow{\pi}\\
    A&\xrightarrow{\quad \mu \quad} &A
   \end{matrix} 
\]
with some map $\psi:S^r\to S^r$ represented by some matrix $M\in S^{r\times r}$. Let $P\in S[T]$ be the characteristic
polynomial of $M$. By the Cayley-Hamilton theorem, we have $P(\psi):S^r\to S^r$ is the zero map. The diagram
\[ \begin{matrix}
    S^r&\xrightarrow{\quad P(\psi)\quad}& S^r\\
    \bigg\downarrow{\pi}& &\bigg\downarrow{\pi}\\
    A&\xrightarrow{\quad P(\mu) \quad} &A
   \end{matrix} 
\]
also commutes. The map $P(\mu):A\to A$ is the multiplication by $P(a)$. Since the $P(\psi)$ is the zero map, the map $P(\mu):A\to A$ is
also zero. Then $P$ is a monic polynomial that annihilates $a$.
\end{proof}

\begin{lemma} \label{lem:mg}
Let $g\in (C\setminus\{0\}$ be a nonzero element of the conductor. Then we have the equality $\langle g\rangle:_S C=g\tilde{S}$.
\end{lemma}

\begin{proof}
$\langle g\rangle:_S C\supseteq g\tilde{S}$: any element in $g\tilde{S}$ can be written as $g\tilde{s}$ for some $\tilde{s}\in\tilde{S}$. If $a\in C$ is 
arbitrary, then $a(g\tilde{s})=(\tilde{s}a)g$ is in $\langle g\rangle$. Hence $g\tilde{S}$ is in $\langle g\rangle:_S C$.

$\langle g\rangle:_S C\subseteq g\tilde{S}$: let $s\in(\langle g\rangle:_S C\setminus\{0\}$. Let $M:=\frac{s}{g}C$. For any $m\in M$, there
is an  $a\in C$ such that $m=a\frac{s}{g}$. For any $\tilde{s}\in\tilde{S}$, we have
\[ m\tilde{s}=a\frac{s}{g}\tilde{s} = (a\tilde{s})\frac{s}{g} = \frac{(a\tilde{s})s}{g} , \]
with $(a\tilde{s})\in C$. Because $s$ is in the quotient ideal of $\langle s\rangle$ by $C$, the numerator is of the form $gs'$ for some $s'\in\tilde{S}$.
Hence $m\tilde{s}\in S$. Since $\tilde{s}$ was chosen arbitrarily in $\tilde{S}$, it follows that $m$ is in the conductor $C$.
So, $\frac{s}{g}C\subseteq C$. By Lemma~\ref{lem:ch}, $\frac{s}{g}\in\tilde{S}$. Therefore $s\in g\tilde{S}$.
\end{proof}

\begin{theorem} Let $g\in C$ of degree $d$, $g\ne 0$. Let $V$ be the vector space of elements of degree $d+1$ in the quotient ideal $I:=\langle g\rangle : C$.
Then the rational map defined by a basis of $V$ maps $B$ to its linear normalization $\tilde{B}$.
\end{theorem}

\begin{proof}
By Lemma~\ref{lem:mg}, $V$ is also the vector space of elements of degree $d+1$ in the ideal $\langle g\rangle_{\tilde{S}}$.
Assume that $(gf_1,\dots,gf_k)$ is a basis of $V$, for some $f_1,\dots,f_k\in\tilde{S}$. Then $f_1,\dots,f_k$ are homogeneous of
degree~1, and they generate the space of all elements of $\tilde{S}$ of degree~1. Therefore the map defined by $(f_1,\dots,f_k)$
maps $B$ to $\tilde{B}$. The map defined by $(gf_1,\dots,gf_k)$ is the same map, because $g$ is may be cancelled in the definition
of a map into projective space.
\end{proof}

Lemma~\ref{lem:mg} and Remark~\ref{rem:nc} imply that Algorithm~\ref{alg:ln} for the computation of the linear normalization of a plane curve is correct.

\floatname{algorithm}{Algorithm}
\renewcommand{\thealgorithm}{LinearNormalization} 
\begin{algorithm}[H]
\caption{}\label{alg:ln}
\begin{algorithmic}[1]
  \Require an irreducible plane curve $B\subset\PP^2$ whose singularities are nodes or cusps.
  \Ensure its linear normalization $\tilde{B}\subset\PP^n$ together with a projection $p:\PP^n\dashrightarrow\PP^3$ mapping
        $\tilde{B}$ to $B$.
  \Statex
  \State {\bf Compute} the radical $I$ of the singular locus of $B$.
  \State {\bf Choose} a homegeneous element $g\in I$ of degree less than $\deg(B)$.
  \State {\bf Compute} the quotient ideal $J:=\langle g,f\rangle:I$, where $f$ is the equation of $B$.
  \State {\bf Extend} $(xg,yg,zg)$ to a basis $(xg,yg,zg,h_1,\dots,h_k)$ for the vector space of forms in $J$ of degree $(\deg(g)+1)$.
  \State {\bf Compute} the image $\tilde{B}\subset\PP^{k+2}$ of $B$ under the rational map $\PP^2\dashrightarrow\PP^{k+2}$ defined by
        this basis.
  \State \Return $\tilde{B}$ together with the projection to the first three coordinates.
\end{algorithmic}
\end{algorithm}

\begin{example}\rm
The nodal cubic curve in Example~\ref{ex:nc} with equation $F=y^2z-x^2z-x^3=0$ has only one node at $(0:0:1)$. Hence we get $I=\langle x,y\rangle$. We choose $g:=x$. The quotient ideal is
\[ J = \langle F,g\rangle:I = \langle x,y^2z\rangle:\langle x,y\rangle = \langle x,yz \rangle.\]
The homogeneous part in degree two is generated by $x^2,xy,xz,yz$. Hence the linear normalization is the image under the map $(x:y:z)\mapsto (x^2:xy:xz:yz)$; 
this is the twisted cubic curve that has been already shown in Example~\ref{ex:nc}.
\end{example}

 \section{Construction and Parametrization of the Veronese surface}

 \subsection{Construction}

Let $B$ be the branching curve of the planar map $\PP^2 \to \PP^2$ defined by polynomials of degree $d\geq 3$. By Proposition~\ref{prop:veronese}, there exists a Veronese surface $V_d$ containing the linear normalization $\tilde{B}$. In this subsection, we will prove that such a Veronese surface is unique, and we give an algorithm for the construction of its ideal. It is known that this ideal is generated by quadratic forms. Note that the Veronese surface $V_d$ is already unique 
up to projective coordinate transformation; but for our purpose it is necessary to compute a Veronese surface containing $\tilde{B}$.

We need to make a case distinction for $d=3$ and $d>3$ -- the latter case is much easier.

 \subsubsection{Degree Three Planar Map}
 
If $d=3$, then $B$ is a curve of degree~18 with 126 double points. The linear normalization $\tilde{B}$ lies in $\PP^9$. Our goal is to find a Veronese surface
$V_3\subset\PP^9$ of degree~9 containing $\tilde{B}$. Proposition~\ref{prop:veronese} tells us a little bit more: the given curve $\tilde{B}$ 
is the image of an irreducible planar sextic under some Veronese map $\nu_3:\PP^2\to \PP^9$. This map is equal to the standard Veronese embedding by cubic monomials 
followed by a projective transformation of $\PP^9$. The Veronese map maps curves of degree $3k$ to intersections with hypersurfaces of degree $k$, for each positive
integer $k$. Hence $\tilde{B}$ is the intersection of $V_3$ with a quadric hypersurface. 

By computation of the ideal of the standard Veronese surface embedded by cubic polynonmials, we see that the ideal of $V_3$ is generated by 27 quadratic forms.
The ideal of the intersection with another quadric hypersurface is generated by 28 quadrics. This is the ideal we have, and what we need is the ideal
generated by the 27 quadratic forms.

It is instructive to compare the Betti tables for the linear normalization and the
Veronese surface of a random example of a planar map defined by degree three polynomials. Here they are:
 \begin{table}[H]
 \centering
 $$\begin{tabular}{c|c c c c c c c c c}
        & 0&1&2&3&4&5&6&7&8 \\ \hline
        total&1&28&132&294&378&294&132&28&1\\
       0: &1&.&.&.&.&.&.&.&. \\
      1:&.&28&105&189&189&105&27&.&.\\
      2:&.&.&27&105&189&189&105&28&.\\
      3:&.&.&.&.&.&.&.&.&1
      \end{tabular}$$
\caption{Betti table for the normalization of the curve.}
\end{table}
\begin{table}[H]
 $$\begin{tabular}{c|c c c c c c c c }
        & 0&1&2&3&4&5&6&7 \\ \hline
        total&1&27&105&189&189&105&27&1\\
       0: &1&.&.&.&.&.&.&. \\
      1:&.&27&105&189&189&105&27&.\\
      2:&.&.&.&.&.&.&.&1
      \end{tabular}$$
\caption{Betti table for the Veronese surface.}
\end{table}
The tables reveal that both ideals have the same linear syzygies (a space of dimension 105). The additional quadratic form in the ideal of $\tilde{B}$
will not show up in any linear syzygy, it is only related by quadratic syzygies to the other 27 quadratic forms. We will make this observation precise
and use it for the uniqueness proof and the algorithm.

Throughout this subsection, we denote by $W_1$ the 28-dimensional vector space of quadratic forms in the ideal of a curve $\tilde{B}\in\PP^9$, which is assumed to be the image of an irreducible sextic plane curve under a Veronese map $v_3$, and we denote by $W_2$ the 27-dimensional subspace of quadratic forms in the ideal of the Veronese surface.

\begin{lemma} \label{lem:thepreviouslystatedlemma}
Let $(Q_1',\dots,Q_{27}')$ be a basis for $W_2$, and let $Q\in W_1\setminus W_2$.
Let $(L_1,L_2,\dots,L_{28})$ be a linear syzygy between quadratic forms $(Q,Q_1',\dots,Q_{27}')$. Then $L_1=0$.
\end{lemma}
 \begin{proof}
Let $f_0,...,f_9$ be the cubic forms defining the Veronese map. We substitute them into the syzygy relation  
$L_1Q+L_2Q_1^{'}+..+L_{28}Q_{27}^{'}=0$.
Since $Q_i^{'}(f_0,..,f_9)=0$ for all $i=1,..,27$, we obtain $L_1(f_0,..,f_9)Q(f_0,..,f_9)=0$. The polynomial $Q(f_0,..,f_9)$ is the pullback of the quadratic form $Q\in I_{\tilde{B}}$, which is the equation of the ramification curve. Furthermore, $Q(f_0,..,f_9)$ is not zero and not a multiple of the equation of the ramification curve since it has degree 6. Therefore, $L_1(f_0,..,f_9)$ must be zero. As a result, $L_1$ vanishes on $V_3$, which implies that $L_1=0$.
\end{proof} 

\begin{lemma} \label{lem:3u}
    The Veronese surface containing $\tilde{B}$ is unique. 
\end{lemma}
\begin{proof}
    Assume, indirectly, that $U_3$ is another Veronese surface containing $\tilde{B}$, $U_3$ not equal $V_3$. 
Let $Q^{''}$ be a quadratic form in the ideal of $U_3$ but not in $W_2$. Since $U_3$ contains $\tilde{B}$, $Q^{''}$ is in $W_1$. By Lemma~\ref{lem:thepreviouslystatedlemma},
any linear syzygy of $Q^{''},Q_1^{'},..,Q_{27}^{'}$ has first entry equal to zero $(*)$. The vectorspace $W_4$ of quadratic forms vanishing on both Veronese surfaces $U_3$ and $V_3$ has dimension 26, by linear algebra. Let $P_1,..,P_{26}$ be a basis. Then  $Q^{''},P_1,..,P_{26}$ is basis
for the vector space of quadratic forms in $I_{U_3}$. By Lemma~\ref{lem:thepreviouslystatedlemma}, the first column of the linear syzygy matrix $S^{'}$ of $(Q^{''},P_1,..,P_{26})^t$
is zero. On the other hand, the $S^{'}$ is the syzygy matrix of a Veronese surface - it is known that the syzygy module of a 
Veronese surface is generated by linear forms. So the first generator $Q^{''}$ is linearly independent from $P_1,..,P_{26}$,
and this absurd. 
\end{proof}

The following algorithm explains each step of constructing the  Veronese surface. 

\floatname{algorithm}{Veronese Surface of}
\renewcommand{\thealgorithm}{Degree3PlanarMap} 
\begin{algorithm}[H]
\caption{}\label{alg:deg3Veronese}
\begin{algorithmic}[1]
  \Require  $\tilde{B}\in\PP^9$, the image of a plane sextic under a Veronese map
  \Ensure the unique Veronese surface $V_3$ containing $\tilde{B}$
  \Statex
  \State {\bf Let} $Q_1,\dots,Q_{28}$ be a set of generators of the ideal of $\tilde{B}$ (they are quadratic)
  \State {\bf Compute} the syzygy matrix S in $k[x_0,..,x_9]^{105\times 28}$ (its entries are linear forms).
  \State {\bf Write } S as $\sum S_i x_i$, $i=0..9$, for some $S_i\in k^{105\times 28}$
  \State {\bf Let} $M$ be the vector space generated by the rows of all $S_i$;
   \State{\bf Compute} a basis for $M$ (27 row vectors in $k^{28}$)
\State {\bf Multiply} each basis vector by column vector $(Q_1,..,Q_{28})^t$ and {\bf return} products.
\end{algorithmic}
\end{algorithm}

\begin{proof}[Proof of correctness]: Let $M$ be the vector space generated by the rows of the matrices $S_i$ as defined in the algorithm. Suppose, indirectly,  that $m\in M$ is such that the scalar product $Q=m\cdot(Q_1,\dots,Q_{28})^t$ is not in the ideal of $V_3$. If $(Q_1^{'},..,Q_{27}^{'})$ is a set of generators for the ideal of $V_3$, then the quadratic forms $(Q,Q_1^{'},..,Q_{27}^{'})$ form another basis for $W_1$. Let $T$ be the transformation matrix from $(Q_1,..,Q_{28})$ to $(Q,Q_1^{'},..,Q_{27}^{'})$ - i.e., $T\cdot (Q_1,..,Q_{28})^t=(Q,Q_1^{'},..,Q_{27}^{'})^t$. Then the syzygy matrix of $(Q,Q_1^{'},..,Q_{27}^{'})$ is $S\cdot T^{-1}$.  Then by Lemma~\ref{lem:thepreviouslystatedlemma}, the first column of the matrix $S\cdot T^{-1}$ is zero. In particular, for each $i\in{0,\dots,9}$, the first column of $S_i\cdot T^{-1}$ is zero. So, for each row of $S_i\cdot T^{-1}$, the first entry is 0. Since $m$ is a linear combination of rows of the matrices $S_i$, we have that $m'=m\cdot T^{-1}$ is a linear combination of rows of $S_i\cdot T^{-1}$. So, the first entry of $m^{'}$ is zero. Therefore, the product $Q=m\cdot(Q_1,..,Q_{28})^t=m\cdot T^{-1}\cdot T\cdot (Q_1,..,Q_{28})^t=m^{'}\cdot (Q,Q_1^{'},..,Q_{27}^{'})^t$ is a linear combination of $Q_1^{'},..,Q_{27}^{'}$, and it follows that $Q$ is in $W_2$ - contradiction.
\end{proof}

\subsubsection{Degree Greater or Equal to Four Planar Map}

If $d>3$, then $\tilde{B}$ is the image of a plane curve $B$ of degree $3(d-1)$ under a Veronese map $v_d:\PP^2\to\PP^N$, where $N=\frac{d(d+3)}{2}$. 

\begin{lemma} \label{lem:4u}
  The Veronese surface $V_d\in\PP^N$ containing $\tilde{B}$ is unique. The ideal of $V_d$ is generated by all quadratic forms in in the ideal of $\tilde{B}$.
\end{lemma}

\begin{proof}
  It suffices to show the second statement. Assume, indirectly, that $\tilde{B}$ is contained in the zero set of a quadratic form $Q$ not vanishing in $V_3$. Then the pullback of $Q$ is a plane curve of degree $2d$ (maybe reducible). This plane curve has to contain the curve $B$ as one of its components. But $2d<3(d-1)$, a contradiction.
\end{proof}

Lemma~\ref{lem:4u} implies that we can compute the ideal of $V_3$ by the following algorithm. 

\floatname{algorithm}{Veronese Surface of}
\renewcommand{\thealgorithm}{DegreeBigger3PlanarMap} 
\begin{algorithm}[H]
\caption{}\label{alg:deg4Veronese}
\begin{algorithmic}[1]
  \Require  $\tilde{B}\in\PP^9$, the image of a plane curve of degree $3(d-1)$ under a Veronese map
  \Ensure the unique Veronese surface $V_d$ containing $\tilde{B}$ 
  \Statex
  \State {\bf Compute} a basis for the vector space of quadratic forms in the ideal of $\tilde{B}$
  \State {\bf Return} this basis
\end{algorithmic}
\end{algorithm}

\subsection{Parametrization of the Veronese surface}

This section focuses on a specific subproblem: given the ideal of a Veronese surface in  $V_d\subset\PP^N$, where $N:=\frac{d(d+3)}{2}$, find a Veronese map $\nu_d:\PP^2\to V_d$  that parametrizes it. In our setting, the Veronese map is defined by an arbitrary basis of the vector space of forms of degree $d$ which is not known, and this makes the subproblem difficult -- if the Veronese map would be defined by monomials, then the subproblem would be trivial.

As mentioned in the introduction, the problem is computationally even harder than linear normalization (see \cite{Schicho:05c}). The main difficulty is a Diophantine subproblem: if we are working over $\QQ$, and we want to find a parametrization with coefficients in $\QQ$, then we need to find a rational point on $V$. In this paper, we do not want to deal with this difficulty. Instead, we assume that we have a black box that can produce points on $V$, and we give a probabilistic algorithm to construct points on $V$ in the case that the ground field is finite. 

Any Veronese map $\nu_d:\PP^2\to V_d$ induces a one-to-one correspondence between linear forms in $\PP^n$ and forms of degree $d$ in $\PP^2$. This will allow as to construct subspaces corresponding to spaces of forms of degree $d$ with useful properties, despite the fact that we do not yet have a Veronese map.

Let $y\in V$ be a point and $k\ge 1$. We define the {\em $k^{th}$ osculating space $O_{y,k}$ to $y$} as the set of all linear forms $L$ such that $L|_{V}$ vanishes to order at least $k$ on $V$ at $y$ (see \cite{BCGI} for a broader treatment of osculating spaces).

\begin{lemma}
Assume that the Veronese map $\nu_d:\PP^2\to V_d$ maps some point $x\in\PP^2$ to $y\in V$. 
Let $k$ be in integer in $\{1,\dots,d\}$. 
Then $O_{y,k}$ corresponds to the vector space of degree $d$ forms that vanish at $x$ to order at least $k$.
\end{lemma}

\begin{proof}
The pullback $\nu_d^{*}$ induces an isomorohism of local rings ${\cal O}_{V,y}$ and ${\cal O}_{\PP^2,x}$, 
and this isomorphism preserves the vanishing order.
\end{proof}

Here is a theorem that allows to construct the inverse of a parametrization of $V$. In order to compute the
a parametrization, we just need to invert the result of that construction.

\begin{theorem} \label{prop:top2}
Let $y_1$ and $y_2$ be two distinct points on the Veronese surface. Let 
\[ J := (O_{y_1,d}\cap O_{y_2,d-1})+(O_{y_1,d-1}\cap O_{y_2,d}) . \]
Then $\dim(J)=3$, and the rational map defined by a basis of $J$ defines an isomorphism $V_d\to\PP^2$.
\end{theorem}

\begin{proof}
Let $\nu_d:\PP^2\to V$ be a Veronese map, and set $x_i:=\nu_d^{-1}(y_i)$ for $i=1,2$. 
Without loss of generality, we assume $x_1=(0:0:1)$ and $x_2=(0:1:0)$.
Then the vectorspace corresponding to $(O_{y_1,d}\cap O_{y_2,d-1})$ is the space of linear forms that
vanish at $x_1$ to order $d$ and at $x_2$ to order $d-1$. This is $\langle s^d,s^{d-1}t\rangle$,
where $s,t,u$ are the homogeneous coordinates of $\PP^2$. Similarily, the vectorspace corresponding to
$(O_{y_1,d-1}\cap O_{y_2,d})$ is $\langle s^d,s^{d-1}u\rangle$. Therefore, the vectorspace corresponding
to $J$ is $\langle s^d,s^{d-1}t,s^{d-1}u\rangle$. The map defined by a basis of $J$ is the composition of
$v_d^{-1}$ with the map $\tau:\PP^2\to\PP^2$ defined by a basis of $\langle s^d,s^{d-1}t,s^{d-1}u\rangle$. 
After cancelling the common factor $s^{d-1}$, we get a basis for the vectorspace of linear forms.
Hence $\tau$ is an isomorphism.
\end{proof}

Theorem~\ref{prop:top2} shows that Algorithm~\ref{alg:par} is correct.

\floatname{algorithm}{Algorithm}
\renewcommand{\thealgorithm}{Parametrization of Veronese Surface}
\begin{algorithm}[H]
\caption{}\label{alg:par}
\begin{algorithmic}[1]
  \Require $V\subset\PP^{\frac{d(d+3)}{2}}$, a Veronese surface
  \Require $y_1,y_2\in V$
  \Ensure a Veronese map $\nu_d:\PP^2\to V_d$
  \Statex
\State {\bf Compute} the osculating spaces $O_{y_1,d-1},O_{y_1,d},O_{y_1,d-1},O_{y_1,d}$
\State {\bf Compute} $J := (O_{y_1,d}\cap O_{y_2,d-1})+(O_{y_1,d-1}\cap O_{y_2,d})$
\State {\bf Let} $\pi:V\to\PP^2$ be the map defined by a basis of $J$
\State \Return $\pi^{-1}:\PP^2\to V$
\end{algorithmic}
\end{algorithm}

If the ground field is finite, then we can find a point on $V$ by the Algorithm~\ref{alg:find}. 

\floatname{algorithm}{Find}
\renewcommand{\thealgorithm}{Points on Surface} 
\begin{algorithm}[H]
\caption{}\label{alg:find}
\begin{algorithmic}[1]
  \Require $V\subset\PP^N$, an algebraic surface
  \Ensure a point $y\in V$
  \Statex
  \State found := {\tt False}
  \While{not found}
   \State {\bf Pick} two random linear forms $L_1,L_2$.
   \State{\bf Compute} the associated primes of $H := I(V)+\langle L_1,L_2\rangle$.
   \If{$H$ has a prime component of a single point $p$}
     \State found := {\tt True}
     \State $y$ := such a point
   \EndIf
  \EndWhile
  \State \Return $y$
\end{algorithmic}
\end{algorithm}

Our computational experiments indicate that the probability of a random pair of linear forms to intersect
$V$ in a point defined over the ground field is $\frac{1}{9}$ for $d=3$ and $\frac{1}{16}$ for $d=4$.
Hence the expected number of loops before finding a point is $9$ for $d=3$ and $16$ for $d=4$.

\begin{remark} \rm
The assumption of the ground field is finite clearly stands against our motivation for the whole reconstruction
problem, namely computer vision. We did the the computations in finite fields mainly to have explicit computational
examples. For computer vision applications, it would actually be a good idea to change to symbolic-numeric
computation, because images have limited resolution leading to numerical errors. Numerically, computing points
on $V$ is not a Diophantine problem and is easier than computing points defined over $\QQ$. However, for many 
other subtasks, such as linear normalization, symbolic-numeric methods are not yet available, and the sensitivity to numerical errors is not known.
\end{remark}

\section{The Whole Algorithm}

For the convenience of the reader, we give here a single algorithm that refers to the functions in the previous section.

\floatname{algorithm}{Algorithm}
\renewcommand{\thealgorithm}{PlanarMap}
\begin{algorithm}[H]
\caption{}\label{alg:planar}
\begin{algorithmic}[1]
  \Require an irreducible plane curve $B\subset\PP^2$ whose singularities are nodes and cusps. Its degree can be factorized as $3d(d-1)$.
  \Ensure a generic planar map $f:\PP^2\to\PP^2$ that has $B$ as a branching curve, if exists.
  \Statex
  \If{$d=2$} use Algorithm ~\ref{alg:case2}.
  \EndIf
  \If{$d\geq 3$}
\State {\bf Compute} the linear normalization $\tilde{B}$ of $B$ by Algorithm~\ref{alg:ln}.
\State{\bf Compute} the Veronese surface containing $\tilde{B}$ for two cases:
      \If{$d=3$}  use the Algorithm  Veronese Surface of ~\ref{alg:deg3Veronese}
\EndIf
      \If{$d\geq 4$} use the Algorithm Veronese Surface of ~\ref{alg:deg4Veronese}) 
     \EndIf 
\State{\bf Find} two points in the Veronese surface $V$; \\ \hfill if ground field is finite, use Algorithm~\ref{alg:find} 
\State{\bf Compute}  parametrization by Algorithm~\ref{alg:par}.
\State{\bf Apply} on Veronese surface V projection in the first three coordinates.
\State \Return the composition of parametrization with the  projection.
\EndIf
\end{algorithmic}
\end{algorithm}

An implementation of degree two can be found in 
\href{https://github.com/HoxhajEriola/Planar-Map-of-Polynomial-of-Degree-Two/tree/main}{\url{https://github.com/HoxhajEriola/Planar-Map-of-Polynomial-of-Degree-Two/tree/main}}

An implementation over finite fields can be found in 
\href{https://github.com/HoxhajEriola/Planar-Map-of-polynomial-degree-3/blob/main/ALG3.m2}{\url{https://github.com/HoxhajEriola/Planar-Map-of-polynomial-degree-3/blob/main/ALG3.m2}} and
\href{https://github.com/HoxhajEriola/Planar-Map-of-Polynomial-of-Degree-Greater-or-Equal-to-4/blob/main/ALG4.m2}{\url{https://github.com/HoxhajEriola/Planar-Map-of-Polynomial-of-Degree-Greater-or-Equal-to-4/blob/main/ALG4.m2}}.

We extensively tested numerous random examples, specifically the field $\mathbb{Z}/32003$. 
We implemented the algorithm in Macaulay2 and tested it on a computer   Intel(R) Core(TM) i5-1135G7 @ 2.40GHz. Below is a table presenting the CPU timings, measured in seconds, for the algorithm steps concerning degree 3 and 4. Algorithm~\ref{alg:case2} was implemented in Maple, it runs in less than 0.3 seconds.

 \begin{table}[H]
  $$\begin{tabular}{|c|c|c|}
  \hline
     Steps   &deg 3&deg 4 \\ \hline
        Adjoint ideal&7.63589
&200.621\\
       Linear Normalization&7.74667&222.676\\
     Veronese Surface &12.927& 218.882 \\
      Points on Veronese&10.9681&225.291\\
      Parametrization&65.0948 & 1145.51\\
      Planar map&87.91753&1126.1\\
      \hline
      \end{tabular}$$
\caption{Timing of the  algorithm for the cases of degree three and four }
\end{table}

\begin{figure}[H]
    \centering

\begin{tikzpicture}
    \node (picture) [] { \includegraphics[height=5cm, width=6cm]{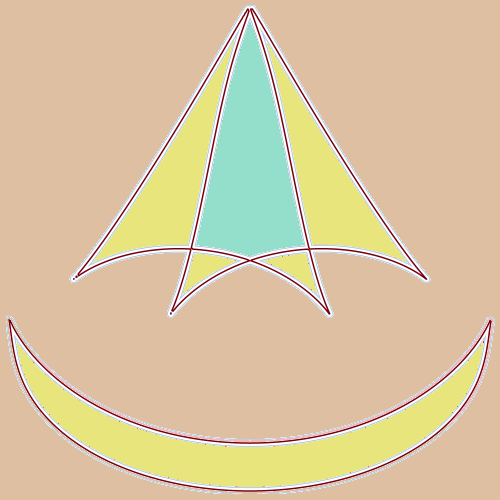}};
    \draw[-latex,red] ($(picture.north west)+(3.1,0)$) arc
    [
        start angle=160,
        end angle=20,
        x radius=2.5cm,
        y radius =1 cm
    ] ;

\end{tikzpicture}
   \includegraphics[height=5cm,width=6cm]{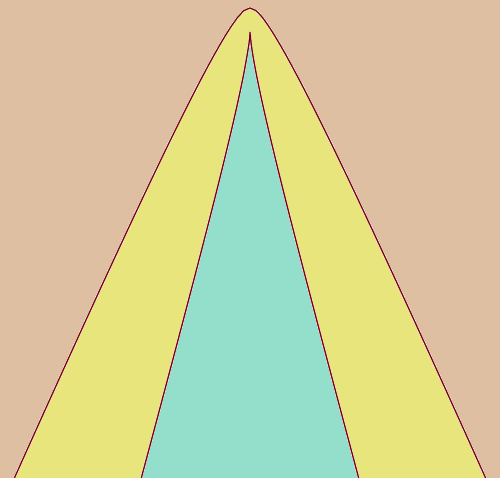} 
    \caption{The branching curve of a planar map of degree~3, with three real nodes and six real cusps. It separates the points on the plane into three distinct regions: the beige area, where each point has three preimage points; the yellow region, where each point has five preimage points; and the blue sky region, where each point has seven preimage points. A point which almost looks like a singularity that is neither a node nor a cusp is magnified -- it is a cusp and a regular part passing close by.}
    \label{fig:deg3_1}
\end{figure}

\begin{figure}[H]
    \centering
   \includegraphics[height=5cm, width=6cm]{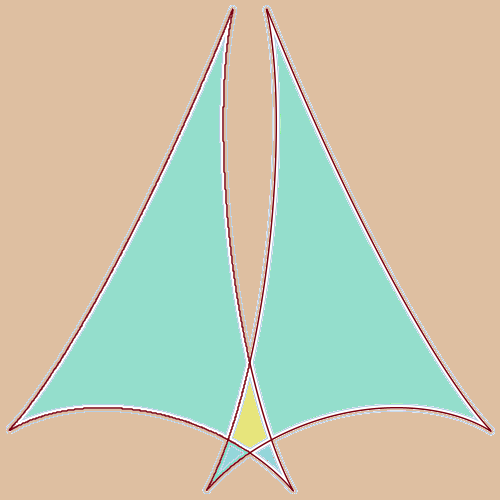} 
 
    \caption{The branching curve of another map of degree~3, with four real nodes and six real cusps. It separates the points on the plane into three distinct regions: the blue sky region, where each point has five preimage points; the yellow region, where each point has six preimage points; and the beige region, where each point has three preimage points.}
    \label{fig:deg3_2}
\end{figure}

The examples~\ref{fig:deg3_1} and \ref{fig:deg3_2} have been computed using the computer algebra system Maple, starting with the map and then computing the branching curve.

\section*{Acknowledgments}

This work is part of GRAPES project that has received funding from the European Union's Horizon 2020 research and innovation programme under the Marie Sklodowska-Curie grant agreement No 860843.
Algorithm~\ref{alg:ln} is based on ideas from a discussion with Matteo Gallet (University of Trieste, Italy). We appreciated the discussions with him and with Niels Lubbes (RICAM, Austrian Academy of Sciences).

\bibliographystyle{amsalpha}

\end{document}